\documentclass[11pt]{amsproc}

\usepackage{eufrak}
\usepackage{amssymb}

\theoremstyle{plain}
\newtheorem{thm}{Theorem}[section]
\newtheorem{lem}[thm]{Lemma}
\newtheorem{prop}[thm]{Proposition}

\theoremstyle{definition}
\newtheorem{rem}[thm]{Remark}
\newtheorem{ntn}[thm]{Notation}
\newtheorem{dft}[thm]{Definition}

\begin{document}

\title[Collections of $p$-subgroups]{New collections of $\mathbf{p}$-subgroups
and homology decompositions for classifying spaces of finite groups}
\maketitle

\begin{center}
{\Large John Maginnis\footnote{Email address: maginnis@math.ksu.edu.} and Silvia Onofrei\footnote{Email address: onofrei@math.ksu.edu.}}\\
{\it Department of Mathematics, Kansas State University,\\ 138 Cardwell
Hall, Manhattan, KS 66506}

\medskip
{\it Communications in Algebra (to appear)}\\
\end{center}

\bigskip
\begin{abstract}
We define new collections of $p$-subgroups for a finite group $G$ and $p$ a prime dividing its order.
We study the homotopy relations among them and with the standard collections of $p$-subgroups and
determine their ampleness and sharpness properties.

\medskip
\noindent{\it Keywords:} Subgroup complexes, homology decompositions, homotopy equivalence.\\
{\it 1991 MSC}: 20J06, 20J05, 51E30, 55P91, 55R35.
\end{abstract}

\bigskip
\section{Introduction}

The paper studies new collections of $p$-subgroups and their relationship to group cohomology and
homology
decompositions, in particular expanding on an earlier systematic study by Grodal and Smith
\cite{grs04}.

\medskip
Let $E_0(G)$ be the set of elements of order $p$ in $G$ which are contained in the center of some
Sylow
$p$-subgroup, and let $E_1(G)$ be the set of elements generated by taking products of commuting
elements in $E_0(G)$.
For $\mathcal{C}$ a collection of $p$-subgroups of $G$, let $\widetilde{\mathcal{C}}$ be the
subcollection
of $\mathcal{C}$ consisting of the subgroups $P \in \mathcal{C}$ such that $Z(P)$ contains a
non-trivial
element in $E_1(G)$, and let $\widehat{\mathcal{C}}$ be the subcollection of
$\widetilde{\mathcal{C}}$
consisting of subgroups $P \in \mathcal{C}$ such that $Z(P)$ contains an element of $E_0(G)$.
We refer to a $p$-subgroup $P \in \widehat{\mathcal{C}}$ as being {\it distinguished}.

\medskip
The main results of this paper are Theorem $3.1$ and Theorem $4.4$,
which state equivariant homotopy equivalences between $\widetilde{\mathcal{A}}_p(G),
\widetilde{\mathcal{B}}_p(G)$ and $\widetilde{\mathcal{S}}_p(G)$, and between
$\widehat{\mathcal{A}}_p(G), \widehat{\mathcal{B}}_p(G)$ and $\widehat{\mathcal{S}}_p(G)$,
where $\mathcal{A}_p$ is the Quillen collection of non-trivial elementary abelian $p$-subgroups,
$\mathcal{B}_p$ is the Bouc collection of non-trivial $p$-radical subgroups, and
$\mathcal{S}_p$ is the Brown collection of non-trivial subgroups.
We also prove homotopy equivalences between certain categories that are related to
homology decompositions.  Theorem $4.4$ contains results on the homotopy
equivalences between these collections and their categories, under one of
three extra assumptions, denoted $(\frak{Ch}), (\frak{Cl})$ and
$(\frak{M})$. In Section 5 we use these equivariant homotopy equivalences,
and apply results of Dwyer \cite{dw98} and Grodal \cite{gr02},
to conclude various sharpness properties of these homology decompositions, which imply that
the mod $p$ homology of the group $G$ fits into an exact sequence, involving the homology of
the subgroups in the collection, or of their centralizers, or of their normalizers.

\medskip
Benson \cite{ben94}
constructed a subcollection $\mathcal{E}_p(G)$ of the Quillen collection.
Benson uses this subcollection for $p=2$ and $G =\text{Co}_3$, in which case it is homotopy
equivalent to a $2$-local geometry $\Delta$ first mentioned by Ronan and
Stroth \cite{rst84}. Although this geometry is not homotopy
equivalent to the Bouc collection, in \cite{mgo} we prove that
$\Delta$ is homotopy equivalent to $\widehat{\mathcal{B}}_2(\text{Co}_3)$,
the collection of $2$-radical subgroups which contain a central
involution lying in the center of some Sylow $2$-subgroup. This result led
us to the definition of distinguished collections of $p$-subgroups.

\medskip
The collection of distinguished $p$-radical subgroups $\widehat{\mathcal{B}}_p(G)$ is contained in
the Bouc collection $\mathcal{B}_p(G)$ and contains a third collection $\mathcal{B}_p^{\rm{cen}}(G)$
which is relevant to homology decompositions and representation theory, and which has been suggested
as a possible generalization of Tits buildings for Lie groups \cite[Sect.1]{y01}. Here
$\mathcal{B}_p^{\rm{cen}}(G)$ is the collection of $p$-centric and $p$-radical subgroups of $G$; see
Section $4.1$ for a precise definition. At the end of Section 4 (in Lemma 4.8) we show that under a
certain (technical) hypothesis, satisfied by all finite Lie groups and several sporadic groups, the
distinguished
Bouc collection equals both the full Bouc collection $\mathcal{B}_p(G)$ and the collection
$\mathcal{B}_p^{\rm{cen}}(G)$. But the example of $\text{Co}_3$ mentioned above shows that these
three collections
need not be homotopy equivalent.  The distinguished collection recovers the homotopy
type of a standard geometry for $\text{Co}_3$ and preserves the points of the
geometry as elements of the collection, the central involutions.

\medskip
In Section $2$ a brief review of basic results from the homotopy theory
of posets and a few facts regarding homology decompositions are given.
In Section $3$ the ``tilde" collections are defined and Theorem $3.1$ is proved. In Section $4$ the
distinguished collections of
$p$-subgroups are discussed. In Section $5$ the sharpness properties of the ``tilde" and
distinguished collections are given.

\medskip
{\bf Acknowledgements:} We would like to thank Stephen Smith for
providing an early version of the manuscript of Benson and Smith \cite{bs04},
the most valuable reference for the material discussed in this paper.
We also thank Masato Sawabe for his careful reading of our final draft.\\

\section{Notations, terminology and standard results}
\subsection{Poset homotopy}
In this section we review a few standard results on the homotopy theory
of posets, most of them due to Quillen \cite{qu78} and proved in the
context of equivariant homotopy by Th\'{e}venaz and Webb \cite{thwb91}.

\medskip
Assume $G$ is a finite group and $p$ a prime dividing its order.
\begin{ntn}Let $\mathcal{X}$ be a poset; set $\mathcal{X}_{>x} =
\lbrace y \in \mathcal{X} \mid y>x \rbrace$. Similarly, we define
$\mathcal{X}_{\geq x}, \mathcal{X}_{< x}, \mathcal{X}_{\leq x}$. Let $|
\mathcal{X} |$  denote the associated simplicial complex with simplices proper inclusion chains in
$\mathcal{X}$. For a $G$-poset
$\mathcal{X}$ and a subgroup $H$ of $G$, $\mathcal{X}^H$ denotes the corresponding fixed point set.
\end{ntn}

\begin{prop}\cite[2.1]{grs04} If $f : \mathcal{X} \rightarrow \mathcal{X}$ is
a $G$-equivariant poset map such that $f \leq Id_{\mathcal{X}}$ (or
dually $f \geq Id_{\mathcal{X}}$), then $\mathcal{X}$ is $G$-homotopy
equivalent to any $G$-invariant subposet $\mathcal{Y}$ which also satisfies
the condition $f(\mathcal{X}) \subseteq \mathcal{Y} \subseteq
\mathcal{X}$.
\end{prop}

\begin{dft}A poset $\mathcal{X}$ is {\it conically contractible} if
there is a poset map $f: \mathcal{X} \rightarrow \mathcal{X}$ and an
element $x_0 \in \mathcal{X}$ such that $x \leq f(x) \geq x_0$ (or dually $x \geq f(x) \leq x_0$) for
all $x \in \mathcal{X}$.
\end{dft}

\begin{lem} A conically contractible poset is contractible. If
$\mathcal{X}$ is conically contractible via a map $f$ which is also
$G$-equivariant then $\mathcal{X}$ is $G$-contractible.
\end{lem}

More generally, if $f,g : \mathcal{X} \rightarrow \mathcal{Y}$ are two poset maps
such that $f(x) \leq g(x)$ for all $x \in \mathcal{X}$, then $f$ and $g$ are homotopic.
This can be thought of as a special case of the homotopy induced by a natural
transformation of functors between categories.

\begin{thm}\cite[1.3]{thwb91} Suppose that $\mathcal{X}$ and $\mathcal{Y}$
are two finite $G$-posets and let $f: \mathcal{X} \rightarrow
\mathcal{Y}$ be a $G$-poset map. The map $f$ is a $G$-homotopy equivalence if and
only if, for all subgroups $H \leq G$, the map $f$ restricts to an
ordinary homotopy equivalence $f^H: \mathcal{X}^H \rightarrow \mathcal{Y}^H$.
\end{thm}

\begin{thm}\cite[Thm. 1]{thwb91} Let $\mathcal{X} \hookrightarrow \mathcal{Y}$ be an inclusion of
$G$-posets. Assume that for all $y \in \mathcal{Y}$ the subposet $\mathcal{X}_{\leq y}$ is
$G_y$-contractible. Then the inclusion is a $G$-homotopy equivalence.
\end{thm}

\subsection{Homology decompositions}
A {\it mod $p$ homology decomposition} of the classifying space $BG$
of a finite group $G$ consists of a mod $p$ homology isomorphism to $BG$ from
a space constructed from classifying spaces of subgroups of $G$, as a homotopy
colimit.  The decomposition is called {\it sharp} if the
associated Bousfield-Kan spectral sequence collapses at the $E^2$-page onto the
vertical axis.

\begin{dft}A {\it collection} $\mathcal{C}$ of subgroups of
$G$ is a set of subgroups which is closed under conjugation. The collection
can be thought of as a poset, or as a category with inclusions as morphisms.
The {\it nerve} or geometric realization $|\mathcal{C}|$ is a simplicial complex.
The group $G$ acts on these structures by conjugation, and the isotropy
group of $H \in \mathcal{C}$ is the normalizer $N_G(H)$.
\end{dft}

A collection is called {\it ample} if the map from the Borel construction
$EG \times_{\substack{G}} |\mathcal{C}| \rightarrow BG$ is a mod $p$
homology isomorphism, where $EG$ is the universal cover of $BG$.

\medskip
Dwyer defines three particular decompositions for an ample collection,
with homotopy colimits that involve the classifying spaces of the
subgroups in the collection, or their centralizers, or their normalizers.
Whether the collection is subgroup sharp, centralizer sharp, or normalizer
sharp can be rephrased in terms of acyclicity for Bredon homology of
the nerves of certain categories, defined below. Details on homology decompositions can be found in
Dwyer \cite{dw98}, Grodal
\cite{gr02} and Benson and Smith \cite[Part I]{bs04}.

\medskip
Let $E\mathbf{A}_{\mathcal{C}}$ denote the
category whose objects are pairs $(H,i)$ where $H$ is a group and $i :
H \rightarrow G$ is a monomorphism with $i(H) \in \mathcal{C}$. A
morphism $(H, i) \rightarrow (K,j)$ is a group homomorphism $\rho: H
\rightarrow K$ with $j \rho = i$. The action of $G$ is given by $g \cdot (H,
i) = (H, c_g i )$, where $c_g: G \rightarrow G$ is given by
$c_g(x)=gxg^{-1}$. The isotropy groups are of the form $C_G(i(H))$.

\medskip
Let $E\mathbf{O}_{\mathcal{C}}$ be the category
whose objects are pairs $(G/H, xH)$ where $H \in \mathcal{C}$ and $xH
\in G/H$. A morphism $(G/H, xH) \rightarrow (G/K, yK)$ is a $G$-map $f:
G/H \rightarrow G/K$ with $f(xH) = yK$. The $G$-action is given by $g
\cdot (G/H, xH) = (G/H, gxH)$ and thus the isotropy groups are of the
form $^xH=xHx^{-1}$. Note that $xHx^{-1} \in \mathcal{C}$.

\medskip
A collection $\mathcal{C}$ of $p$-subgroups is {\it centralizer
sharp, subgroup sharp}, or {\it normalizer sharp} if the nerves of $E
\mathbf{A}_{\mathcal{C}}, E \mathbf{O}_{\mathcal{C}}$ or $\mathcal{C}$ are
acyclic for Bredon homology. A comprehensive analysis of the ampleness and sharpness properties
of nine collections of $p$-subgroups can be found in Grodal and Smith \cite{grs04}.

\begin{rem} The three categories $E\mathbf{O}_{\mathcal{C}}, \;
\mathcal{C}$ and $E\mathbf{A}_{\mathcal{C}}$ are equivalent and therefore they
have homotopy equivalent nerves. In general these equivalences are not
$G$-homotopy equivalences; see \cite[Section 3.8]{dw98}. For a subgroup $H$ of $G$, the fixed point
sets of the above categories are related by the equivalences
$E\mathbf{O}_{\mathcal{C}}^H \rightarrow \mathcal{C}_{\geq H}$ and
$E\mathbf{A}_{\mathcal{C}}^H \rightarrow \mathcal{C}_{ \leq C_G(H)}$; see \cite[Section 2]{grs04}.
\end{rem}

\begin{lem}[\cite{grs04}, $2.5$] Let $G$ be a finite group and
$\mathcal{C}' \subset \mathcal{C}$ collections of subgroups of $G$.
\begin{enumerate}
\item[(1)]Assume for all $P \in \mathcal{C} \setminus \mathcal{C}'$
that $\mathcal{C}_{>P}$ is contractible. Then
$|E\mathbf{O}_{\mathcal{C}'}| \rightarrow |E\mathbf{O}_{\mathcal{C}}|$ is a $G$-homotopy
equivalence.
\item[(2)]Assume for all $P \in \mathcal{C} \setminus \mathcal{C}'$
that $\mathcal{C}_{<P}$ is contractible. Then
$|E\mathbf{A}_{\mathcal{C}'}| \rightarrow |E\mathbf{A}_{\mathcal{C}}|$ is a $G$-homotopy
equivalence.
\item[(3)]Assume that $\mathcal{C}_{>P}$ is
$N_G(P)$-contractible for all $P \in \mathcal{C} \setminus \mathcal{C}'$. Then $|\mathcal{C}'|
\rightarrow |\mathcal{C}|$ is a
$G$-homotopy equivalence.
\end{enumerate}
\end{lem}

\bigskip
\section{Collections related to the Benson collection}

In the sequel, several collections of $p$-subgroups are defined; they generalize the collection
$\mathcal{E}_p(G)$ introduced by Benson \cite{ben94} in order to study the mod-$2$ cohomology of the
group Co$_3$. These ``tilde" collections consist of $p$-subgroups whose centers contain elements from
the set $E_1(G)$ defined below. The main result of this section is contained in Theorem $3.1$ and
says that the corresponding ``tilde" complexes have homotopy properties which are similar to those of
their standard counterparts, the Quillen, Bouc and Brown complexes.

\medskip
Consider the smallest subset $E_1(G)$ which
satisfies the following properties:
\begin{enumerate}
\item[\phantom{ii}i)]contains the elements of order $p$ in the center
of a Sylow $p$-subgroup of $G$;
\item[\phantom{i}ii)]is closed under conjugation in $G$;
\item[iii)]is closed under taking products of commuting elements.
\end{enumerate}

Then $\mathcal{E}_p(G)$ is the collection of elementary abelian
$p$-subgroups of $G$ which are subsets of $E_1(G)$. We will call
$|\mathcal{E}_p(G)|$ the Benson complex.

\medskip
Let $\mathcal{C}_p (G)$ be any collection of $p$-subgroups of $G$. For
a $p$-subgroup $P$ of $G$ define:
$$\widetilde{P} = \Omega _1 Z(P) \cap E_1(G)$$
which is a group since $E_1(G)$ is closed under commuting
products. Denote:
$$\widetilde{\mathcal{C}}_p(G) = \lbrace P | P \in \mathcal{C}_p(G) \;
\text{and} \; \widetilde{P} \not = 1 \rbrace \subseteq \mathcal{C}_p
(G)$$

The following collections are standard in the literature:
\begin{align*}
\mathcal{A}_p (G) &= \lbrace E \mid E \; \text{nontrivial elementary
abelian p-subgroup of} \; G  \rbrace;\\
\mathcal{S}_p (G) &= \lbrace P \mid P \; \text{nontrivial p-subgroup
of} \; G \rbrace;\\
\mathcal{B}_p (G) &= \lbrace R \mid R \; \text{nontrivial p-radical
subgroup of}\; G \rbrace.
\end{align*}

The simplicial complex $| \mathcal{A}_p(G)|$ is known as the Quillen
complex, $| \mathcal{S}_p(G) |$ is known as the Brown complex and
$|\mathcal{B}_p(G)|$ is called the Bouc complex. The three complexes are
$G$-homotopy equivalent \cite[Thm.2]{thwb91}. A $p$-subgroup
$R$ of $G$ is called $p$-{\it radical} if $R$ is the largest normal
$p$-subgroup of its normalizer, that is $R = O_p(N_G(R))$.

\begin{thm} Let $\mathcal{C}$ be one of the collections
$\mathcal{E}_p(G), \; \widetilde{\mathcal{A}}_p(G), \; \widetilde{\mathcal{S}}_p(G)$
or $\widetilde{\mathcal{B}}_p(G)$. Then there exist certain homotopy equivalences, summarized in the
following table:\\

\begin{picture}(1000, 50)(0,0)
\put(118,25){$\mathcal{E}_p(G)$}
\put(188,25){$\widetilde{\mathcal{A}}_p(G)$}
\put(258,25){$\widetilde{\mathcal{S}}_p(G)$}
\put(328,25){$\widetilde{\mathcal{B}}_p(G)$}

\put(65,-10){$|E \mathbf{O}_{\mathcal{C}}|$}
\put(122,-10){$\bullet$}
\put(129,-10){$\cdots$}
\put(144,-10){$\cdots$}
\put(159,-10){$\cdots$}
\put(174,-10){$\cdots$}
\put(192,-10){$\bullet$}
\put(199,-10){$\cdots$}
\put(214,-10){$\cdots$}
\put(229,-10){$\cdots$}
\put(244,-10){$\cdots$}
\put(262,-10){$\bullet$}
\put(263,-7){\line(1,0){70}}
\put(332,-10){$\bullet$}

\put(65,-55){$|\mathcal{C}|$}
\put(122,-55){$\bullet$}
\put(123,-47){$\vdots$}
\put(123,-35){$\vdots$}
\put(123,-23){$\vdots$}
\put(123,-52){\line(1,0){70}}
\put(192,-55){$\bullet$}
\put(193,-47){$\vdots$}
\put(193,-35){$\vdots$}
\put(193,-23){$\vdots$}
\put(193,-52){\line(1,0){70}}
\put(262,-55){$\bullet$}
\put(263,-42){$|$}
\put(263,-24){$|$}
\put(263,-52){\line(1,0){70}}
\put(332,-55){$\bullet$}
\put(333,-42){$|$}
\put(333,-24){$|$}

\put(65,-100){$|E \mathbf{A}_{\mathcal{C}}|$}
\put(122,-100){$\bullet$}
\put(123,-87){$|$}
\put(123,-69){$|$}
\put(123,-97){\line(1,0){70}}
\put(192,-100){$\bullet$}
\put(193,-97){\line(1,0){70}}
\put(193,-87){$|$}
\put(193,-69){$|$}
\put(262,-100){$\bullet$}
\put(263,-87){$|$}
\put(263,-69){$|$}
\put(269,-100){$\cdots$}
\put(284,-100){$\cdots$}
\put(299,-100){$\cdots$}
\put(314,-100){$\cdots$}
\put(332,-100){$\bullet$}
\put(333,-92){$\vdots$}
\put(333,-80){$\vdots$}
\put(333,-68){$\vdots$}

\put(210,-120){{\it Table} $\mathit{3.1}$}
\end{picture}
\end{thm}

\vspace*{4cm}
\begin{ntn} A solid line corresponds to a $G$-homotopy equivalence, a
dashed line to an $S$-homotopy equivalence and a dotted line to an
ordinary homotopy equivalence. Here $S$ denotes a Sylow $p$-subgroup of $G$. The categories $E
\mathbf{A}_{\mathcal{C}}$ and $E \mathbf{O}_{\mathcal{C}}$ were defined in the previous section.
\end{ntn}

\begin{proof}[{\bf Proof of Theorem 3.1}]The proof of the Theorem consists of a number of steps, at
each step (except for $\mathit{(3)}$ which is a technical result) we prove one (or more) of the
homotopy equivalences corresponding to the horizontal and vertical lines in {\it Table}
$\mathit{3.1}$. The solid horizontal lines from the
center row are given by steps $\mathit{(1)}$, $\mathit{(2)}$ and $\mathit{(4)}$ below.

\medskip
{\it $\mathit{(1)}$. The inclusion $\mathcal{E}_p (G) \hookrightarrow
\widetilde{\mathcal{A}}_p (G)$ is a $G$-homotopy equivalence}.\\
Let $P \in \widetilde{\mathcal{A}}_p (G)$, then:
$$\mathcal{E}_p (G) _{\leq P} = \lbrace E \in \mathcal{E}_p (G)| E \leq
P \rbrace = \widetilde{\mathcal{A}}_p (\widetilde{P}) = \mathcal{A}_p
(\widetilde{P})$$
is $N_G(P)$-contractible since it is a cone on $\widetilde{P}$. The conclusion follows by an
application of Theorem $2.6$.

\medskip
{\it $\mathit{(2)}$. The inclusion $\widetilde{\mathcal{A}}_p
(G)\hookrightarrow \widetilde{\mathcal{S}}_p (G)$ is a $G$-homotopy equivalence}.\\
We show that $\widetilde{\mathcal{A}}_p(G)_{\leq P}$ is
contractible for any $P \in \widetilde{\mathcal{S}}_p (G)$ and we apply
Theorem $2.6$ again. If $Q \in \widetilde{\mathcal{A}}_p (G)_{\leq P}$
then $Q \widetilde{P} \in \widetilde{\mathcal{A}}_p(G)_{\leq P}$ and the double inequality $Q \leq Q
\widetilde{P} \geq
\widetilde{P}$ provides the conical contraction (see Lemma $2.4$).
The poset map $f(Q)=Q \widetilde{P}$ is $N_G(P)$-equivariant.

\medskip
{\it $\mathit{(3)}$. Let $P \in \widetilde{\mathcal{S}}_p(G)$ and $Q \in
\mathcal{S}_p(G)$ satisfy $P \trianglelefteq Q$. Then $Q \in
\widetilde{\mathcal{S}}_p(G)$}.\\
Since the subgroup $\widetilde{P}$ is characteristic in $P$,
it follows that $\widetilde{P} \trianglelefteq Q$ and thus $\widetilde{P} \cap Z(Q) \not =
1$; hence $\widetilde{Q} \not = 1$.

\medskip
{\it $\mathit{(4)}$. The inclusion $\widetilde{\mathcal{B}}_p(G) \hookrightarrow
\widetilde{\mathcal{S}}_p(G)$ is a $G$-homotopy equivalence}.\\
Assume that $P \in \widetilde{\mathcal{S}}_p(G) \setminus
\widetilde{\mathcal{B}}_p(G)$ and $Q \in
\widetilde{\mathcal{S}}_p(G)_{>P}$. Then $P < N_Q(P) \leq N_G(P)$ and also $P < O_p (N_G(P)) \leq
N_G(P)$. According to step $\mathit{(3)}$, the subgroups $N_Q(P), O_p(N_G(P))$ and $N_Q(P)
O_p(N_G(P))$ are all in $\widetilde{\mathcal{S}}_p(G)$. Hence we
obtain the chain:
$$Q \geq N_Q(P) \leq N_Q(P) O_p(N_G(P)) \geq O_p(N_G(P))$$
of poset maps $\widetilde{\mathcal{S}}_p(G)_{>P} \rightarrow
\widetilde{\mathcal{S}}_p(G)_{>P}$. The poset maps $f_1(Q)=N_Q(P)$ and $f_2(R)=R \cdot O_p(N_G(P))$
satisfy $f_1(Q) \leq Q$ and $f_2(R) \geq R$.  As these
poset maps are $N_G(P)$-equivariant, it follows that $\widetilde{\mathcal{S}}_p(G)_{>P}$ is
$N_G(P)$-contractible and an application of Lemma $2.9(3)$ gives the result.

\medskip
$\mathit{(5)}$. For the left solid horizontal line on the bottom row, consider $P \in
\widetilde{\mathcal{A}}_p(G) \setminus \mathcal{E}_p(G)$ and note that
$\widetilde{\mathcal{A}}_p(G)_{<P}$ is conically contractible via $Q \geq
\widetilde{Q} \leq \widetilde{P}$ where $P> Q \in
\widetilde{\mathcal{A}}_p(G)$. Then apply the result of Lemma $2.9(2)$.

\medskip
$\mathit{(6)}$. Given $P \in \widetilde{\mathcal{S}}_p(G) \setminus
\widetilde{\mathcal{B}}_p(G)$, the poset $\widetilde{\mathcal{S}}_p(G)_{>P}$ is
contractible; see the proof of step $\mathit{(4)}$. Thus, according to Lemma $2.9(1)$, $|E
\mathbf{O}_{{\widetilde{\mathcal{B}}}_p}|$ and $|E \mathbf{O}_{{\widetilde{\mathcal{S}}}_p}|$ are
$G$-homotopy equivalent, and this takes care of the solid line on the top row.

\medskip
$\mathit{(7)}$. Next we prove the $G$-homotopy equivalence corresponding to the middle
solid horizontal line on the bottom row. Set $\mathcal{C} =
\widetilde{\mathcal{S}}_p(G)$ and $\mathcal{C}' =
\widetilde{\mathcal{A}}_p(G)$. Let $H$ be a subgroup of $G$. Applying Theorem $2.5$
to $E\mathbf{A}_{\mathcal{C}'} \rightarrow  E \mathbf{A}_{\mathcal{C}}$, we show that
$E\mathbf{A}_{\mathcal{C}'}^H \rightarrow  E \mathbf{A}_{\mathcal{C}}^H$ is a homotopy
equivalence. Recall from Remark $2.8$ that there is an
equivalence of categories $E \mathbf{A}^H_{\mathcal{C}} \rightarrow
\mathcal{C}_{\leq C_G(H)}$; thus it suffices to show that $\mathcal{C}_{\leq C_G(H)}$ and
$\mathcal{C}'_{\leq C_G(H)}$ are homotopy equivalent. To simplify the notation we let $\mathcal{X} =
\mathcal{C}'_{\leq C_G(H)}$ and $\mathcal{Y} = \mathcal{C}_{\leq C_G(H)}$. The proof uses a
nonequivariant
version of Theorem $2.6$; we show that $\mathcal{X}_ {\leq P}$ is
contractible for any $P \in \mathcal{Y}$. First note $\widetilde{P} \in
\mathcal{X}$. If $Q \in \mathcal{X}_{\leq P}$ then $Q \widetilde{P} \in
\mathcal{X}_{\leq P}$ also. The conical contractibility $\mathcal{X}_{\leq P}$ follows from the
double
inequality $Q\leq  Q \widetilde{P} \geq \widetilde{P}$ which corresponds to
poset maps.  Thus the inclusion $\mathcal{X} \hookrightarrow \mathcal{Y}$ is
a homotopy equivalence.

\medskip
$\mathit{(8)}$. The dashed vertical line in the $\mathcal{E}_p(G)$ column, follows from
\cite[Theorem 1.1]{grs04}. The lower dashed vertical lines from the columns
corresponding to $\widetilde{\mathcal{A}}_p(G)$ and
$\widetilde{\mathcal{S}}_p(G)$ are obtained by combining the parallel vertical dashed lines
with the adjacent solid horizontal lines.

\medskip
{\it $\mathit{(9)}$. Let $S$ be a fixed Sylow $p$-subgroup of $G$. For $\mathcal{C} =
\widetilde{\mathcal{S}}_p(G)$, there is an $S$-homotopy equivalence between
$|E{\mathbf{O}}_{\mathcal{C}}|$ and $| \mathcal{C}|.$}\\
Applying Theorem $2.5$, we first prove that for any subgroup $H$ of $S$, the
functor $E{\mathbf{O}}_{\mathcal{C}}^H \rightarrow \mathcal{C}^H$ is
a homotopy equivalence. As $E\mathbf{O}_{\mathcal{C}}^H$ is equivalent to $\mathcal{C}_{\geq H}$ (see
Remark $2.8$), we have to show that the inclusion $\mathcal{C}_{\geq H}
\rightarrow {\mathcal{C}}^H$ is a homotopy equivalence. Let
$Q \in \mathcal{C}^H$; thus $Q \in \widetilde{\mathcal{S}}_p(G)$ and $H
\leq N_G(Q)$. Hence $Q \leq QH \leq N_G(Q)$, and by step $\mathit{(3)}$ we have $QH \in
\widetilde{\mathcal{S}}_p(G)$. Also $QH \in
\mathcal{C}^H$ since $H$ normalizes $QH$. It follows that the map $Q
\rightarrow QH$ is a poset endomorphism on $\mathcal{C}^H$ with image in
$\mathcal{C}_{\geq H}$ and by Proposition $2.2$, $\mathcal{C}^H$ and $\mathcal{C}_{\geq H}$ are
homotopy equivalent. The top dashed line from column
$\widetilde{\mathcal{S}}_p(G)$ follows now.

\medskip
$\mathit{(10)}$. The dashed line in the
$\widetilde{\mathcal{B}}_p(G)$ column is obtained by composing the solid horizontal lines and
the dashed line from the upper right corner rectangle.
\end{proof}

\medskip
\begin{rem}We give an example which shows that the three dotted
horizontal lines cannot be replaced by dashed (or solid) lines. Let $G= D_8$.
This group contains a unique central involution and four involutions of
noncentral type. There are two subgroups of the form $\mathbb{Z}/2
\times \mathbb{Z}/2$ and a subgroup $\mathbb{Z}/4$. Therefore
$\mathcal{E}_2(G)$ contains a unique point and $\widetilde{\mathcal A}_2(G)$ is formed of two
edges sharing a common vertex. Let $H $ be one of the rank two
elementary abelian subgroups, then $\mathcal{E}_2(G)_{\geq H}$ is empty while
$\widetilde{\mathcal{A}}_2(G)_{\geq H}$ is a point. This shows the
dotted line on the left of the top row in {\it Table} $\mathit{3.1}$ cannot be replaced by a
dashed line. Next, note that $C_{G}(H) = H$ and thus
$\widetilde{\mathcal{S}}_2(G)_{\leq C_G(H)}$ is an edge, but
$\widetilde{\mathcal{B}}_2(G)_{\leq C_G(H)}$ is empty. This takes care of the dotted line in the
right bottom row.  Finally, let $K=\mathbb{Z}/4$. Then
$\widetilde{\mathcal{A}}_2(G)_{\geq K}$ is empty, while
$\widetilde{\mathcal{S}}_2(G)_{\geq K}$ contains an edge, thus it is contractible. This explains the
middle dotted line on the top row.
\end{rem}

\section{Distinguished collections of $p$-subgroups}

In this section additional new collections are defined; they consist of $p$-subgroups which contain
in their centers elements which also lie in the center of some Sylow $p$-subgroup of $G$. We analyze
the homotopy relations between these ``distinguished" collections as well as between their
corresponding
categories $E{\mathbf O}_{\mathcal{C}}$ and $E {\mathbf A}_{\mathcal{C}}$. It turns out that these
relations are quite similar to those among the standard collections if $G$ satisfies one of three
conditions formulated below. These results are contained in Theorem $4.4$.

\medskip
Assume that $G$ is a finite group and $p$ a prime dividing its order.
Denote by $E_0(G)$ the family of elements of order $p$ in $G$ of central type, that is those which
lie in the center of some Sylow $p$-subgroup of $G$.

\medskip
For a $p$-subgroup $P$ of $G$ define:
$$\widehat{P} = \langle x | x \in \Omega _1 Z(P) \cap E_0(G)
\rangle$$
Further, for $\mathcal{C}_p(G)$ a collection of $p$-subgroups of $G$
denote:
$$\widehat{\mathcal{C}}_p(G) = \lbrace P | P \in \mathcal{C}_p(G) \;
\text{and} \; \widehat{P} \not= 1\; \rbrace$$
We call $\widehat{\mathcal{C}}_p(G)$ the {\it distinguished}
$\mathcal{C}_p(G)$ {\it collection}. We shall refer to the subgroups in
$\widehat{\mathcal{C}}_p(G)$ as {\it distinguished subgroups}.

\begin{rem} When $G$ has one conjugacy class of elements of order $p$, it
is obvious that $\widehat{\mathcal{C}}_p(G) = \mathcal{C}_p(G)$. If $P \leq G$ is a $p$-subgroup then
$\widehat{P} \leq
\widetilde{P}$; also $\widehat{\mathcal{C}}_p(G) \subseteq
\widetilde{\mathcal{C}}_p(G)$.
\end{rem}

Before proceeding to the study of the distinguished collections of
$p$-subgroups we formulate the following conditions:
\begin{enumerate}
\item[($\mathfrak{M}$)] \it{Given any $P \in
\widehat{\mathcal{S}}_p(G)$, the subgroup $N_G(P)$ is contained in a $p$-local subgroup which
contains a Sylow $p$-subgroup of $G$.}
\item[($\mathfrak{Cl}$)] \it{The elements of order $p$ and of central type in $G$
are closed under taking products of commuting elements.}
\item[($\mathfrak{Ch}$)]\it{The group $G$ is of local characteristic
$p$-type; this means that for every $p$-local subgroup $H$ of $G$ the following holds: $C_H(O_p(H))
\leq O_p(H)$.}
\end{enumerate}

\begin{rem}Note that if $R$ is a $p$-subgroup of $G$ and if $N_G(R)$
contains a Sylow $p$-subgroup of $G$, then $R$ is distinguished. This is
easy to see since in this case $R \triangleleft S \in
\text{Syl}_p(N_G(R)) \subseteq \text{Syl}_p(G)$. Thus $R \cap Z(S) \not= 1$.
\end{rem}

\begin{rem}[{\it The three conditions and the sporadic groups}]The
three conditions listed above were formulated with the properties of the
$26$ sporadic simple groups in mind. For $p=2$ all but four of the sporadic simple groups satisfy at
least one of the three assumptions. The sporadic
groups of local characteristic $2$-type are: $\text{M}_{11}$, $\text{M}_{22}$,
$\text{M}_{23}$, $\text{M}_{24}$, $\text{J}_1$, $\text{J}_3$,
$\text{J}_4$, $\text{Co}_2$ and $\text{Th}$. The sporadic groups which
satisfy $(\frak{M})$ are $\text{M}_{12}$, $\text{M}_{24}$, $\text{J}_1$,
$\text{J}_2$, $\text{J}_4$, $\text{McL}$, $\text{Suz}$, $\text{Co}_3$,
$\text{Co}_1$, $\text{Th}$, $\text{M}$, $\text{He}$ and $\text{O'N}$.
The condition $(\mathfrak{Cl})$ holds for the following groups: $\text{M}_{12}$, $\text{J}_2$,
$\text{HS}$, $\text{Suz}$,
$\text{Co}_3$ and $\text{Ru}$; each of these six groups has two classes of involutions. The sporadic
groups $\text{M}_{11}$, $\text{M}_{22}$, $\text{M}_{23}$, $\text{J}_1$,
$\text{J}_3$, McL, Th, O'N, Ly have one class of involutions each and the assumption
$(\mathfrak{Cl})$ clearly holds; in all these cases the distinguished collections are the same as the
standard ones. The group $\text{Fi}_{23}$ has three classes of involutions, all central, and
therefore satisfies the $(\mathfrak{Cl})$ assumption, too. The sporadic groups
which do not satisfy any of the three conditions are $\text{Fi}_{22}$,
$\text{Fi}'_{24}$, HN and BM.
\end{rem}

\begin{thm} Let $\mathcal{C}$ be one of the collections
$\widehat{\mathcal{A}}_p(G), \; \widehat{\mathcal{S}}_p(G)$ or
$\widehat{\mathcal{B}}_p(G)$. There exist
homotopy equivalences, summarized in the following table:

\begin{picture}(1000, 50)(0,0)
\put(140,25){$\widehat{\mathcal{A}}_p(G)$}
\put(230,25){$\widehat{\mathcal{S}}_p(G)$}
\put(320,25){$\widehat{\mathcal{B}}_p(G)$}

\put(65,-10){$|E \mathbf{O}_{\mathcal{C}}|$}
\put(142,-10){$\bullet$}
\put(149,-10){$\cdots$}
\put(164,-10){$\cdots$}
\put(179,-10){$\cdots$}
\put(194,-10){$\cdots$}
\put(209,-10){$\cdots$}
\put(224,-10){$\cdots$}
\put(232,-10){$\bullet$}
\put(233,-7){\line(1,0){90}}
\put(253,-7){$\phantom{a}^{(\mathfrak{Cl}, \;\mathfrak{Ch},
\;\mathfrak{M})}$}
\put(322,-10){$\bullet$}

\put(65,-55){$|\mathcal{C}|$}
\put(142,-55){$\bullet$}
\put(143,-52){\line(1,0){90}}
\put(143,-47){$\vdots$}
\put(143,-35){$\vdots$}
\put(143,-23){$\vdots$}
\put(232,-55){$\bullet$}
\put(233,-42){$|$}
\put(233,-24){$|$}
\put(233,-35){$\phantom{a}^{(\mathfrak{Cl}, \;\mathfrak{Ch})}$}
\put(233,-52){\line(1,0){90}}
\put(253,-52){$\phantom{a}^{(\mathfrak{Cl}, \;\mathfrak{Ch}, \;
\mathfrak{M})}$}
\put(322,-55){$\bullet$}
\put(323,-42){$|$}
\put(323,-24){$|$}
\put(323,-35){$\phantom{a}^{(\mathfrak{Cl}, \;\mathfrak{Ch})}$}

\put(65,-100){$|E \mathbf{A}_{\mathcal{C}}|$}
\put(142,-100){$\bullet$}
\put(143,-87){$|$}
\put(143,-69){$|$}
\put(143,-80){$\phantom{a}^{(\mathfrak{Cl}, \;\mathfrak{Ch})}$}
\put(142,-97){\line(1,0){90}}
\put(232,-100){$\bullet$}
\put(233,-87){$|$}
\put(233,-69){$|$}
\put(233,-80){$\phantom{a}^{(\mathfrak{Cl}, \;\mathfrak{Ch})}$}
\put(239,-100){$\cdots$}
\put(254,-100){$\cdots$}
\put(269,-100){$\cdots$}
\put(284,-100){$\cdots$}
\put(299,-100){$\cdots$}
\put(314,-100){$\cdots$}
\put(322,-100){$\bullet$}
\put(323,-92){$\vdots$}
\put(323,-80){$\vdots$}
\put(323,-68){$\vdots$}

\put(210,-120){{\it Table} $\mathit{4.4}$}
\end{picture}
\end{thm}

\vspace*{4cm}

\begin{ntn} The same notation as in  {\it Table} $\mathit{3.1}$ is used. A label $(\mathfrak{c})$
means that the corresponding homotopy equivalence holds under hypothesis $(\mathfrak{c})$. The
notation for the standard collections of $p$-subgroups is given at the beginning of Section $3$ and
the categories $E \mathbf{A}_{\mathcal{C}}$ and $E \mathbf{O}_{\mathcal{C}}$ were defined in Section
$2$.
\end{ntn}

\begin{proof}[{\bf Proof of Theorem 4.4}] The proof of the Theorem consists of a number of steps,
given below. The homotopy relations corresponding to the two solid horizontal lines from the center
row are proved in steps $\mathit{(1)}$ and $\mathit{(3)}$.

\medskip
{\it $\mathit{(1)}$. The inclusion $\widehat{\mathcal{A}}_p(G) \hookrightarrow
\widehat{\mathcal{S}}_p(G)$ is a $G$-homotopy equivalence.}\\
For $P \in \widehat{\mathcal{S}}_p(G)$ the subposet
$\widehat{\mathcal{A}}_p (G)_{\leq P}$ is $N_G(P)$-contractible via the
double inequality $Q \leq Q \widehat{P} \geq \widehat{P}$ given by
an $N_G(P)$-equivariant poset map, where $Q \in \widehat{\mathcal{A}}_p(G)_{\leq P}$. The
$G$-homotopy equivalence follows by an
application of Theorem $2.6$.

\medskip
{\it $\mathit{(2)}$. Let $P \in \mathcal{S}_p(G)$ and $Q \in
\widehat{\mathcal{S}}_p(G)_{>P}$. Then $N_Q(P)\in \widehat{\mathcal{S}}_p(G)$.}\\
For $Q \in \widehat{\mathcal{S}}_p(G)$ with $P<Q$, $P < N_Q(P)\leq Q$ and $Z(Q) \leq Z(N_Q(P))$.
Hence $N_Q(P) \in \widehat{\mathcal{S}}_p(G)$.

\medskip
{\it $\mathit{(3)}$. Assume that $G$ satisfies one of the
conditions $(\frak{Cl})$, $(\frak{Ch})$ or $(\frak{M})$. Then $\widehat{\mathcal{S}}_p(G)$ and
$\widehat{\mathcal{B}}_p(G)$ are $G$-homotopy equivalent.}\\
Under each of the three assumptions on $G$ and for a subgroup $P \in \widehat{\mathcal{S}}_p(G)
\setminus \widehat{\mathcal{B}}_p(G)$ we show that $\widehat{\mathcal{S}}_p(G)_{>P}$ is
$N_G(P)$-contractible. Then we use Lemma $2.9(3)$ to obtain the result.

\smallskip
Denote $O_{NP} = O_p(N_G(P))$; since $P$ is not $p$-radical it follows that $P < O_{NP}$. Let $Q \in
\widehat{\mathcal{S}}_p(G)_{>P}$. For a Sylow $p$-subgroup $\bar{S}$ of $N_G(P)$ which contains
$N_Q(P)$, let $S$ denote a Sylow
$p$-subgroup of $G$ such that $\bar{S} = S \cap N_G(P)$. Since $Z(S)
\leq C_G(P) \leq N_G(P)$ it follows that $Z(S) \leq Z(\bar{S})$. Also $O_{NP} \trianglelefteq
\bar{S}$.

\smallskip
$(\mathfrak{M})$. Let $M$ be a $p$-local subgroup of $G$ with $N_G(P) \leq
M$ and such that $M$ contains a Sylow $p$-subgroup of $G$.  Denote $R = O_p(M)$ and assume that $S$,
the Sylow $p$-subgroup chosen in the previous paragraph, lies in $M$.
Since $R \trianglelefteq S$ the intersection $Z(S) \cap R$ is
nontrivial. Consider the following string of inequalities:
$$Q \geq N_Q(P)  \leq N_Q(P) N_R(P) O_{NP} \geq N_R(P) O_{NP}$$
The subgroup $N_Q(P)$ normalizes $N_R(P)$ since $R \trianglelefteq M$ and therefore $N_R(P)
\trianglelefteq N_G(P)$. Hence $N_Q(P)N_R(P)$ is a $p$-subgroup of $N_G(P)$ which normalizes
$O_{NP}$, since $O_{NP} \trianglelefteq N_G(P)$. From the fact that $1 \not= Z(S) \cap R \leq
Z(N_R(P)
O_{NP})$ it follows that $N_R(P) O_{NP}$ is distinguished. Also $1 \not =
Z(S) \cap R \leq Z(N_Q(P)N_R(P) O_{NP})$ and hence this last
$p$-subgroup is also distinguished. Since all the inequalities
correspond to $N_G(P)$-equivariant poset maps $\widehat{\mathcal{S}}_p(G)_{>P} \rightarrow
\widehat{\mathcal{S}}_p(G)_{>P}$ it follows that $\widehat{\mathcal{S}}_p(G)_{>P}$ is
$N_G(P)$-contractible.

\smallskip
$(\mathfrak{Ch})$. Since $G$ is of local characteristic $p$-type
$C_{N_G(P)}(O_{NP}) \leq O_{NP}$. Thus $P< O_{NP} \leq \bar{S} \leq S$ and $Z(S)
\leq C_{N_G(P)}(O_{NP}) \leq O_{NP}$. Hence $Z(S) \leq Z(O_{NP})$ and $O_{NP} \in
\widehat{\mathcal{S}}_p(G)$. It also follows $N_Q(P)O_{NP}$ is distinguished, since $Z(S) \leq
Z(N_Q(P)O_{NP})$. Now consider the string of poset maps $\widehat{\mathcal{S}}_p(G)_{>P} \rightarrow
\widehat{\mathcal{S}}_p(G)_{>P}$
$$Q \geq N_Q(P) \leq N_Q(P) O_{NP} \geq O_{NP}$$
which proves the $N_G(P)$-contractibility of $\widehat{\mathcal{S}}_p(G)_{>P}$ in this case.

\smallskip
$(\mathfrak{Cl})$. Under this condition, $\widehat{\mathcal{S}}_p(G) = \widetilde{\mathcal{S}}_p(G)$.
The result follows from step $\mathit{(4)}$ of Theorem $3.1$.

\medskip
$\mathit{(4)}$. For $\mathcal{C}' = \widehat{\mathcal{A}}_p(G)$ and $\mathcal{C} =
\widehat{\mathcal{S}}_p(G)$, the proof that $|E\mathbf{A}_{\mathcal{C}'}|$ and
$|E\mathbf{A}_{\mathcal{C}}|$ are $G$-homotopy equivalent is similar to that in step $\mathit{(7)}$
of Theorem $3.1$. Let
$H$ be a subgroup of $G$ and denote $\mathcal{X} = \mathcal{C}'_{\leq
C_G(H)}$ and $\mathcal{Y} =\mathcal{C}_{\leq C_G(H)}$. We show that
$\mathcal{X}_{\leq P}$ is $C_G(P)$-contractible for any $P \in \mathcal{Y}$.
For $Q \in \mathcal{X}_{\leq P}$ we obtain the contracting homotopy $Q
\leq \widehat{P} Q \geq \widehat{P}$ which proves the conical
contractibility of $\mathcal{X}_{\leq P}$. Now apply Theorem $2.6$ to show that $\mathcal{X}$ and
$\mathcal{Y}$ are homotopy equivalent, and Remark $2.8$ and Theorem $2.5$ to obtain the
$G$-homotopy equivalence between $|E\mathbf{A}_{\mathcal{C}}|$ and
$|E\mathbf{A}_{\mathcal{C}'}|$.

\medskip
$\mathit{(5)}$. Next consider $\mathcal{C}' = \widehat{\mathcal{B}}_p(G)$ and
$\mathcal{C} = \widehat{\mathcal{S}}_p(G)$. Let $P \in \mathcal{C} \setminus
\mathcal{C}'$. From step $\mathit{(3)}$ we know that if $G$
satisfies one of the hypotheses $(\frak{Cl})$, $(\frak{Ch})$ or $(\frak{M})$,
the subposet  $\mathcal{C}_{>P}$ is $N_G(P)$-contractible. Next we
apply Lemma $2.9(1)$ to obtain the $G$-homotopy equivalence between
$|E\mathbf{O}_{\mathcal{C}'}|$ and $|E\mathbf{O}_{\mathcal{C}}|$.

\medskip
$\mathit{(6)}$. If we assume
$(\frak{Cl})$ the two vertical dashed lines from the $\widehat{\mathcal{S}}_p(G)$
column follows from Theorem $3.1$ and from the fact that
$\widetilde{\mathcal{S}}_p(G) = \widehat{\mathcal{S}}_p(G)$ in this case.

\medskip
$\mathit{(7)}$. {\it Assume that} $(\frak{Ch})$ {\it holds in $G$. Then $\mathcal{B}_p(G) =
\widehat{\mathcal{B}}_p(G)$ and the collections $\mathcal{S}_p(G)$ and $\widehat{\mathcal{S}}_p(G)$
are $G$-homotopy equivalent.}\\
Assume that $G$ is of local characteristic $p$-type. Since $\widehat{\mathcal{B}}(G) \subseteq
\mathcal{B}_p(G)$, it remains to prove that if $R$ is a $p$-radical subgroup in $G$, then $R$ is also
a distinguished $p$-subgroup. The center $Z(S)$ of any Sylow $p$-subgroup $S$ of $G$ which contains
$R$ lies in $C_G(R)$. Since $G$ satisfies $(\frak{Ch})$ and $R$ is $p$-radical,
$C_{N_G(R)}(O_p(N_G(R))=C_G(R) \leq R$ and therefore $Z(S) \leq C_G(R) = Z(R)$, proving that $R$ is
distinguished.

\smallskip
In step $\mathit{(3)}$ it was proved that if $G$ is of local characteristic $p$-type, the collections
$\widehat{\mathcal{B}}_p(G)$ and $\widehat{\mathcal{S}}_p(G)$ are $G$-homotopy equivalent. Also
$\mathcal{B}_p(G)$ and $\mathcal{S}_p(G)$ are $G$-homotopy equivalent; see \cite[Thm.2]{thwb91}. Thus
$\mathcal{S}_p(G)$ and $\widehat{\mathcal{S}}_p(G)$ are also $G$-homotopy equivalent in this case.

\medskip
$\mathit{(8)}$. Assume that $(\frak{Ch})$ holds. Then $\widehat{\mathcal{B}}_p(G) = \mathcal{B}_p(G)$
and the upper
vertical dashed line in the column $\widehat{\mathcal{B}}_p(G)$ follows from
\cite[Theorem 1.1]{grs04};  the upper dashed line in the
$\widehat{\mathcal{S}}_p(G)$ column follows now by
composing the corresponding existing adjacent lines.

\medskip
$\mathit{(9)}$. Next consider the lower dashed vertical line in the
$\widehat{\mathcal{S}}_p(G)$ column under the assumption $(\frak{Ch})$. Let $S$ be a Sylow
$p$-subgroup of $G$. Denote $\mathcal{C} = \mathcal{S}_p(G)$ and
$\mathcal{C}'=\widehat{\mathcal{S}}_p(G)$. Recall that $|\mathcal{C}|$ and
$|E\mathbf{A}_{\mathcal{C}}|$ are $S$-homotopy equivalent
\cite[Theorem 1.1]{grs04} and that the inclusion $\mathcal{C}' \hookrightarrow
\mathcal{C}$ induces a $G$-homotopy equivalence on nerves, according to step $\mathit{(7)}$. In order
to show that $|\mathcal{C}'| \rightarrow |E
\mathbf{A}_{\mathcal{C}'}|$ is an $S$-homotopy equivalence, it suffices
to prove that $|E\mathbf{A}_{\mathcal{C}}|$ and
$|E\mathbf{A}_{\mathcal{C}'}|$ are $S$-homotopy equivalent. We use Theorem $2.5$ and Remark $2.8$.
For $P \leq S$, the subposet
$\mathcal{C}_{\leq C_G(P)}$ is contractible, since $O_p(C_G(P))>1$. We show
that $\mathcal{C}'_{\leq C_G(P)}$ is also contractible. The subgroup $O_{NP} = O_p(N_G(P))$ is
distinguished; see step $\mathit{(3)}$ above. Hence $Z(O_{NP}) \in \mathcal{C}'_{\leq C_G(P)}$. For
$Q \in
\mathcal{C}'_{\leq C_G(P)}$ the contracting homotopy $Q \leq Q Z(O_{NP}) \geq
Z(O_{NP})$ proves the contractibility of $\mathcal{C}'_{\leq C_G(P)}$.
Thus the inclusion $\mathcal{C}'_{\leq C_G(P)} \rightarrow
\mathcal{C}_{\leq C_G(P)}$ induces a homotopy equivalence. It follows that
$|E\mathbf{A}_{\mathcal{C}}|$ and $|E\mathbf{A}_{\mathcal{C}'}|$ are
$S$-homotopy equivalent and $|\mathcal{C}'|$ and $|E\mathbf{A}_{\mathcal{C}'}|$
are $S$-homotopy equivalent, too.

\medskip
$\mathit{(10)}$. The dashed vertical line in column $\widehat{\mathcal{A}}_p(G)$, under one of the
assumptions $(\frak{Ch})$ or $(\frak{Cl})$ follows now from composing the horizontal
solid lines with the vertical dashed line in the lower left rectangle of
the table.
\end{proof}

\begin{rem} To see that the two remaining dotted horizontal lines cannot be
replaced by dashed lines, consider the example $G=D_8$ from Remark $3.3$. This group satisfies all
three hypotheses
$(\frak{M}), (\frak{Cl})$ and $(\frak{Ch})$. The ``tilde'' collections are the
same as the distinguished collections, since the group satisfies $(\frak{Cl})$
and all the facts mentioned therein are applicable in this case, too. For the dotted upper vertical
line in the left column, let $\mathcal{C} =
\widehat{\mathcal{A}}_2(D_8)$ and $H=\mathbb{Z}_4$. Then
$E\mathbf{O}_{\mathcal{C}}^H \simeq \mathcal{C}_{\geq H}$ is empty but $\mathcal{C}^H$ is
contractible. Let
now $\mathcal{C} = \widehat{\mathcal{B}}_2(D_8)$ and $H =
\mathbb{Z}_4$. In this case $E\mathbf{A}_{\mathcal{C}}^H \simeq \mathcal{C}_{\leq C_G(H)}$ is empty
and $\mathcal{C}^H$
is contractible. Thus the dotted vertical line in the right
column cannot be replaced by a dashed line.
\end{rem}

\subsection{Collections of $\mathbf{p}$-centric subgroups} In the last part of this section, we
describe some relations between the distinguished collections and collections of $p$-centric
subgroups. The
$p$-subgroup $R$ is $p$-{\it centric} if $Z(R)$ is a Sylow $p$-subgroup
of $C_G(R)$. Let $\mathcal{C}e_p(G)$ denote the subcollection of $\mathcal{S}_p(G)$
consisting of $p$-centric subgroups and let
$\mathcal{B}^{\text{cen}}_p(G)=\mathcal{C}e_p(G) \cap \mathcal{B}_p(G)$ be the collection of
nontrivial $p$-radical and $p$-centric subgroups. These two collections are
not in general homotopy equivalent with $\mathcal{S}_p(G)$, but the inclusion map
$\mathcal{B}^{\text{cen}}_p(G) \subseteq \mathcal{C}e_ p(G)$
induces a $G$-homotopy equivalence.

\medskip
The collection $\mathcal{B}_p^{\rm{cen}}(G)$ was studied by various authors; see for example Dwyer
\cite[Sect.10]{dw98} and Grodal and Smith \cite{grs04}. The relations between
$\mathcal{B}_2^{\text{cen}}(G)$ and certain known $2$-local geometries were analyzed in Benson and
Smith \cite[Chp.8]{bs04}, and Sawabe \cite{sa03} constructs a $p$-local geometry based on the
collection of $p$-centric and $p$-radical subgroups. Therefore, it seems natural to investigate the
relation between the distinguished Bouc collection and the collection of $p$-radical and $p$-centric
subgroups.

\begin{prop} a). All the $p$-centric subgroups of $G$ are distinguished.\\
b). The collection of distinguished $p$-radical subgroups
contains the collection of $p$-centric and $p$-radical subgroups:
$\mathcal{B}_p^{\text{cen}}(G) \subseteq \widehat{\mathcal{B}}_p(G)$.
\end{prop}

\begin{proof} a). Let $P \in \mathcal{C}e_p(G)$ and let $S$ be any Sylow
$p$-subgroup of $G$ which contains $P$. Then $Z(S) \leq C_G(P)$, hence $Z(S) \leq Z(P)$ and $P \in
\widehat{\mathcal{S}}_p(G)$.\\
b). Follows from part a).
\end{proof}

\begin{lem} If $G$ is of local characteristic $p$-type then
$\mathcal{B}_p(G) = \widehat{\mathcal{B}}_p(G) =
\mathcal{B}^{\text{cen}}_p(G)$.
\end{lem}

\begin{proof} The first equality was proved in step $\mathit{(7)}$ of Theorem $4.4$. It remains to
show that every distinguished $p$-radical subgroup of $G$ is also $p$-centric. Let $R$ be a
$p$-radical subgroup of $G$
thus $R = O_p(N_G(R))$. Since $G$ is of local characteristic $p$-type, $R$
contains the group
$$C_{N_G(R)}(O_p(N_G(R))) =C_{N_G(R)}(R) = C_G(R) \cap N_G(R) =
C_G(R).$$
It follows now that $C_G(R) = Z(R)$ proving that $R$ is $p$-centric.
\end{proof}

\begin{rem}In general, the inclusions $\mathcal{B}^{\text{cen}}_p(G) \subseteq
\widehat{\mathcal{B}}_p(G) \subseteq \mathcal{B}_p(G)$ do not induce homotopy equivalences. To see
this consider the example of the sporadic simple group
$G=$Co$_3$ and $p=2$. This group satisfies both conditions $(\mathfrak{Cl})$
and $(\mathfrak{M})$. We proved in \cite{mgo} that
$\widehat{\mathcal{B}}_2(G)$ is not homotopy equivalent to $\mathcal{B}_2(G)$ or
to $\mathcal{B}_2^{\text{cen}}(G)$.
\end{rem}

\begin{rem} Recall from Remark $4.3$, that the sporadic groups
which do not satisfy any of the three assumptions are $\text{Fi}_{22}$,
$\text{Fi}'_{24}$, HN and BM. However, in all these four cases
$\widehat{\mathcal{B}}_2(G) = \mathcal{B}_2^{\text{cen}}(G)$.
\end{rem}

\section{Ampleness and sharpness properties}
The following {\it Table} $\mathit{5.1}$ summarizes the sharpness results for the
distinguished collections defined in Section 4, consisting of
groups whose centers contain a nontrivial element in the center
of a Sylow $p$-subgroup of $G$, and for the
``tilde" collections defined in Section 3, consisting of groups
whose centers contain a nontrivial element of the set $E_1(G)$
(generated by products of commuting elements of central type). Definitions of ampleness and sharpness
were given in Section $2$.\\

\begin{center}
\begin{tabular}{|c||c|c|c||c|c|c|}
\hline
\phantom{\Huge{A}}&$\quad \widetilde{\mathcal{A}}_p \quad $&$\quad
\widetilde{\mathcal{B}}_p \quad $&$\quad  \widetilde{\mathcal{S}}_p \quad
$&$\quad  \widehat{\mathcal{A}}_p \quad $&$\quad
\widehat{\mathcal{B}}_p \quad $&$\quad  \widehat{\mathcal{S}}_p \quad $\\
\hline
s&n&y&y&n&y&y\\
\hline
n&y&y&y&y&y&y\\
\hline
c&y&n&y&y&n&y\\
\hline
\end{tabular}

\vspace*{.5cm}
{\it Table 5.1}
\end{center}

\begin{thm}The ``tilde" and the distinguished collections are all ample and have the following sharpness properties:
\begin{itemize}
\item[(a)] $\widetilde{\mathcal{A}}_p(G)$ and
$\widehat{\mathcal{A}}_p(G)$ are centralizer and normalizer sharp.
\item[(b)] $\widetilde{\mathcal{B}}_p(G)$ and
$\widehat{\mathcal{B}}_p(G)$ are subgroup and normalizer sharp.
\item[(c)] $\widetilde{\mathcal{S}}_p(G)$ and
$\widehat{\mathcal{S}}_p(G)$ are subgroup, centralizer and normalizer sharp.
\end{itemize}
\end{thm}

\begin{proof}[{\bf Proof of Theorem 5.1}]
Grodal and Smith \cite{grs04} show that equivariant homotopy equivalences
can be used to propagate sharpness properties via Bredon cohomology.  A
$G$-homotopy equivalence induces an isomorphism in Bredon cohomology.
An $S$-homotopy equivalence, for $S$ a Sylow $p$-subgroup of $G$, is
sufficient to propagate sharpness since the composition of the restriction
map and the transfer map equals multiplication by the index, which is
relatively prime to $p$, and so the composition is an isomorphism; see \cite[Thm.6.4]{dw98}.

\medskip
The Benson collection $\mathcal{E}_p(G)$, those nontrivial elementary abelian
$p$-subgroups $P$ of $G$ satisfying $P \subseteq E_1(G)$, is centralizer sharp \cite{ben94}. Grodal
and Smith use an $S$-homotopy
equivalence to conclude it is also normalizer sharp; see \cite[Thm. 1.2]{grs04}.  This same
technique
applies in our situation to show that
the $G$ and $S$-homotopy equivalences in
Theorem $3.1$ imply all of the sharpness properties for
$\widetilde{\mathcal{A}}_p(G) , \widetilde{\mathcal{B}}_p(G)$ and
$\widetilde{\mathcal{S}}_p(G)$. For example, the $G$-homotopy equivalences
between all four of these collections (the solid lines in the middle row
of {\it Table} $\mathit{3.1}$) imply that the three ``tilde" collections are normalizer
sharp.  Then certain $S$-homotopy equivalences (the vertical dashed lines of
{\it Table} $\mathit{3.1}$) yield the stated subgroup or centralizer sharpness properties.

\medskip
The fact that the distinguished collections $\widehat{\mathcal{B}}_p(G)$ and
$\widehat{\mathcal{S}}_p(G)$ are both subgroup and normalizer sharp follows
immediately from work of Grodal.
In Grodal \cite{gr02}, the
sharpness properties of the collection $\mathcal{D}_p(G)$ of principal
$p$-radical subgroups are investigated; these are $p$-centric subgroups $P$ (see Section $4.1$) with
the property that $O_p(N_G(P)/PC_G(P))=1$. As these groups are by definition both
$p$-centric and $p$-radical, $\mathcal{D}_p(G) \subseteq
\widehat{\mathcal{B}}_p(G) \subseteq
\widehat{\mathcal{S}}_p(G)$ by Proposition $4.7$. Then Grodal's Theorems $1.2$ and $7.3$
from \cite{gr02} imply that $\widehat{\mathcal{S}}_p(G)$ and
$\widehat{\mathcal{B}}_p(G)$ are subgroup and normalizer sharp. Further,
$\widehat{\mathcal{A}}_p(G)$ is normalizer sharp since, by Step $\mathit{(1)}$ in the proof of
Theorem $4.4$, it is
$G$-homotopy equivalent to $\widehat{\mathcal{S}}_p(G)$.

\medskip
Section $8$ of the paper of Dwyer \cite{dw98} discusses the sharpness
properties of elementary abelian $p$-subgroups, and his Theorem $8.3$
reduces the question of the centralizer sharpness to the subcollection of
groups contained in a fixed Sylow $p$-subgroup $S$. So we consider the collection
$\mathcal{C}= \widehat{\mathcal{A}}_p(G)_{\leq S}$ of distinguished elementary
abelian $p$-subgroups of $G$ which lie in $S$. To verify the centralizer sharpness
of $\widehat{\mathcal{A}}_p(G)$, it is sufficient to show that the nerve of the
category $E \mathbf{A}_{\mathcal{C}}$ is contractible.
The conical contraction of $E \mathbf{A}_{\mathcal{C}}$ is obtained using the distinguished
elementary abelian subgroup $Z = \Omega _1 Z(S)$. Denote the
inclusion map by $j:Z \hookrightarrow S$. Given any monomorphism $i:H
\rightarrow S$, with image $i(H)$ a distinguished elementary abelian p-group
(the pair $(H,i)$ is an object in $E \mathbf{A}_{\mathcal{C}}$), construct the subgroup $H' = i(H)
\cdot Z$, with the
corresponding inclusion map $i': H' \rightarrow S$. Note that $H'$ is elementary
abelian, and is distinguished since $Z \leq H'$.
Then we have a zigzag of natural transformations:
$$(H,i) \rightarrow (H', i') \leftarrow (Z, j)$$
proving the contractibility of $E \mathbf{A}_{\mathcal{C}}$ and the centralizer sharpness of
$\widehat{\mathcal{A}}_p(G)$.

\medskip
The centralizer sharpness of $\widehat{\mathcal{S}}_p(G)$ now follows from
the $G$-homotopy equivalence between $E \mathbf{A}_{\mathcal{C}'}$ and
$E\mathbf{A}_{\mathcal{C}}$ proven in Step $\mathit{(4)}$ of Theorem $4.4$ for $\mathcal{C}'
= \widehat{\mathcal{A}}_p(G)$ and $\mathcal{C} =
\widehat{\mathcal{S}}_p(G)$.
\end{proof}

\begin{rem}
The quaternion group $Q_8$ of order $8$ has a periodic mod
$2$ cohomology, which is not detected on the central
$\mathbb{Z}_2$ (there are nilpotent cohomology classes in the kernel of the
restriction map). This implies that for the group $Q_8$, the collection
$\widehat{\mathcal{A}}_2(Q_8)$ is not subgroup sharp, and the collection
$\widehat{\mathcal{B}}_2(Q_8)$ consisting only of the group $Q_8$ is
not centralizer sharp.
\end{rem}

\end{document}